\documentclass[12pt]{article}
\usepackage[psamsfonts]{amsfonts}
\usepackage{amsmath,amssymb,amsthm}
\usepackage[dvips]{graphicx}
\usepackage{eepic}
\newtheorem{THM}{Theorem}[section]
\newtheorem{COR}[THM]{Corollary}

\newtheorem{LEM}[THM]{Lemma}
\newtheorem{PRP}[THM]{Proposition}

\newtheorem{CON}[THM]{Conjecture}

\def\UseSection{
        \numberwithin{equation}{section}
	\theoremstyle{plain}
        \newtheorem{theorem}    {Theorem}[section]
        \DefineTheorems 
}

\def\DefineTheorems{
	
	\newtheorem{lemma}      [theorem] {Lemma}
	
	\newtheorem{prop}       [theorem] {Proposition}
	
	\newtheorem{cor}        [theorem] {Corollary}

	\theoremstyle{definition}
	\newtheorem{defn}       [theorem] {Definition}

	\theoremstyle{definition}

}

\newcommand{\bt}   {\begin{theorem}}
\newcommand{\et}   {\end  {theorem}}
\newcommand{\bl}   {\begin{lemma}}
\newcommand{\el}   {\end  {lemma}}
\newcommand{\bp}   {\begin{prop}}
\newcommand{\ep}   {\end  {prop}}
\newcommand{\bc}   {\begin{cor}}
\newcommand{\ec}   {\end  {cor}}
\newcommand{\bd}   {\begin{defn}}
\newcommand{\ed}   {\end  {defn}}

\newcommand{\ba}   {\begin{array}}
\newcommand{\ea}   {\end  {array}}
\newcommand{\be}   {\begin{enumerate}}
\newcommand{\ee}   {\end  {enumerate}}
\newcommand{\bi}   {\begin{itemize}}
\newcommand{\ei}   {\end  {itemize}}

\def\eq#1\en{\begin{equation}#1\end{equation}}  
\def\eqsplit#1\ensplit{
	\begin{equation}\begin{split}#1\end{split}\end{equation}
	}
\def\eqalign#1\enalign{
	\begin{align}#1\end{align}
	}
\def\eqmul#1\enmul{
	\begin{multline}#1\end{multline}
	}
\newcommand{\eqarrstar} {\begin{eqnarray*}} 
\newcommand{\enarrstar} {\end{eqnarray*}} 
\newcommand{\eqarray}   {\begin{eqnarray}} 
\newcommand{\enarray}   {\end{eqnarray}}

\newcommand{\lbeq}[1]  {\label{e:#1}}
\newcommand{\refeq}[1] {\eqref{e:#1}}    

\makeatletter
\newcommand{\labelcounter}[2]{{%
	\stepcounter{#1}
	\protected@write\@auxout{}%
	{\string\newlabel{#2}{{\csname the#1\endcsname}{\thepage}}}%
	{\ref{#2}}
	}}
\makeatother
\newcommand{\sss}   { \scriptscriptstyle } 

\newcommand{\Qbold} {{\mathbb Q}}
\newcommand{\Rbold} {{\mathbb R}}

\newcommand{\Zbold} {{\mathbb Z}}

\newcommand{\Zd}    {{ {\Zbold}^d }}

\newcommand{\spose}[1] {{\hbox to 0pt{#1\hss}} }
\newcommand{\ltapprox} {\mathrel{\spose{\lower 3pt\hbox{$\mathchar"218$}}
 \raise 2.0pt\hbox{$\mathchar"13C$}}}
\newcommand{\gtapprox} {\mathrel{\spose{\lower 3pt\hbox{$\mathchar"218$}}
 \raise 2.0pt\hbox{$\mathchar"13E$}}}

\newcommand{\prodtwo}[2]{
	\prod_{ \mbox{ \scriptsize 
		$\begin{array}{c} 
		{#1} \\ 
		{#2}  
		\end{array} $ } 
		} 
	} 

\UseSection  
\usepackage[usenames]{color}

\newcommand{\ch}[1]{#1}
\newcommand{\ci}[1]{#1}
\newcommand{\blank}[1]{ }
\newcommand{\ra}{\rightarrow}

\newcommand{\bigcupd} {\stackrel{\cdot}{\bigcup}}
\newcommand{\del}{\partial}

\newcounter{countC}  
\newcounter{countR}  
\newcommand{\R}{\Rbold}
\newcommand{\Z}{\Zbold}

\newcommand{\nn}{\nonumber}

\newcommand{\smallsup}[1] {{\scriptscriptstyle{({#1}})}}

\newcommand{\walk}{\vec{\omega}}
\newcommand{\ewalk}{\vec{\eta}}
\newcommand{\walkvec}[2]{\vec{\omega}^{\smallsup{#1}}_{#2}}
\newcommand{\walkcoor}[2]{\omega^{\smallsup{#1}}_{#2}}

\newcommand{\wh}[1]{\widehat{#1}}
\newcommand{\mc}[1]{\mathcal{#1}}
\newcommand{\mP}{\mathbb{P}}
\newcommand{\mQ}{\mathbb{Q}}
\newcommand{\mE}{\mathbb{E}}

\newcommand{\eqn}[1]{\eq #1 \en}
\newcommand{\sN}{{\sss N}}

\title  {Excited against the tide:\\
  A random walk with competing drifts
        }

\author{
Mark Holmes\footnote{Department of Statistics,
The University of Auckland, Private Bag 92019, Auckland 1142,
New Zealand. E-mail {\tt mholmes@stat.auckland.ac.nz}}
}

\begin{document}

\maketitle

    \begin{abstract}
    We study a random walk that has a drift $\frac{\beta}{d}$ to the right when located at a previously unvisited vertex and a drift $\frac{\mu}{d}$ to the left otherwise.  We prove that in high dimensions, for every $\mu$, the drift to the right is a strictly increasing and continuous function of $\beta$, and that there is precisely one value $\beta_0(\mu,d)$ for which the resulting speed is zero.
    \end{abstract}


\section{Introduction}
\label{sec-intro}
In this paper we study what might be called an {\it excited random walk with drift} (ERWD),  where the random walker has a drift $\frac{\beta}{d}$ 
in the positive direction of the first component each time the walker visits a
new site, and a drift $\frac{\mu}{d}$ in the negative direction of the first component on subsequent visits.  Models of this type with drift for a fixed finite number of visits to a site have been studied by Zerner and others, usually in 1 dimension, see for example \cite{AR05,BS08,BS08b,Zern05}.  They are generalisations of the {\em excited random walk} introduced in \cite{BW03}. 

It is known that excited random walk has ballistic behaviour when $d\geq 2$ in
\cite{BW03,Kozm03, Kozm05}, while there is no ballistic behaviour (when $\beta<1$) in one dimension \cite{Dav99}.  Laws of large numbers and central limit theorems can be obtained for $d\ge2$ using renewal techniques (see for example \cite{SZ99}, \cite{Zern05}, \cite{BR07}).
Intuitively, the velocity appearing in the laws of large numbers should be increasing in the excitement parameter, $\beta$.  This has been proved in dimensions $d\ge 9$ \cite{HH08mono} using a perturbative expansion developed in \cite{HH07}, but the problem remains open for $d<9$.  Using the same expansion approach it is possible \cite{H08rwre} to prove monotonicity for the first coordinate of the speed of random walk in a partially random i.i.d.~environment, in the special case where at each vertex, either the left or the right step is not available.

We use the same argument in this paper, together with the law of large numbers provided by \cite[Theorem 1.4]{BSZ03} to prove that for $d\ge 12$, for each fixed $\mu>0$ the speed in the direction of the positive first coordinate is continuous and strictly increasing in $\beta$.  We give an easy coupling argument showing that the speed is negative when $\beta$ is sufficiently small (depending on $\mu$ and $d$) and show that when $d\ge 9$ the speed is positive when $\beta$ is sufficiently large.  We conclude that when $d\ge 12$, for each $\mu\ge 0$ there is exactly one value $\beta_0(\mu,d)$ of $\beta$ for which the speed is zero.

Although we only consider the {\em once}-excited random walk with all subsequent visits having a reverse drift, many of the results of this paper will remain valid under appropriate relaxations of these conditions, e.g.~multi-excited random walks with a reverse drift only after the $k$th visit to a site, or random walks in site-percolation-like cookie environments.  It may also be possible (using the approach in \cite{HH08mono,H08rwre} and this paper) to prove monotonicity in $p$ for the speed of excited random walk in a site-percolation (parameter $p$) cookie environment in high dimensions.

\subsection{Main results and organisation}
\label{sec-results}
A nearest-neighbour
random walk path $\vec{\eta}$ is a sequence
$\{\eta_i\}_{i=0}^\infty$ for which $\eta_i\in \Z^d$
and $\eta_{i+1}-\eta_i$ is a nearest-neighbour of the origin
for all $i\geq 0$.
For a general nearest-neighbour path $\vec{\eta}$ with $\eta_0=0$, we write
$p^{\vec{\eta}_i}(x_i, x_{i+1})$
for the conditional probability that the walk steps from $\eta_i=x_i$ to
$x_{i+1}$, given the history of the path $\vec{\eta}_i=(\eta_0, \ldots, \eta_i)$.
We write $\walk_n$ for the $n$-step path of excited random walk with opposing drift (ERWD)
and $\Qbold$ for the law of $\{\walk_n\}_{n=0}^{\infty}$, i.e.,
for every $n$-step nearest-neighbour path $\ewalk_n$,
    \eq
    \lbeq{SIRP}
    \Qbold(\walk_n=\vec{\eta}_n)
    =\prod_{i=0}^{n-1} p^{\vec{\eta}_i}(\eta_{i},\eta_{i+1}),
    \en
where, for $i=0$, $p^0(0,\eta_{1})=(2d)^{-1}(1+\beta e_1\cdot \eta_1)$ is the probability that the first step is to 
$\eta_1$, and
\eq
    p^{\vec{\eta}_i}(\eta_{i},\eta_{i}+x)=\frac{1+e_1\cdot x(\beta I_{\{\eta_{i}\notin \vec{\eta}_{i-1}\}}-\mu I_{\{\eta_{i}\in \vec{\eta}_{i-1}\}})}{2d}=\frac{1+e_1\cdot x((\beta+\mu) I_{\{\eta_{i}\notin \vec{\eta}_{i-1}\}}-\mu))}{2d}.
\en
Here $e_1=(1,0, \ldots, 0)$, and $x\cdot y$ is the
inner-product between $x$ and $y$.  It is always the case that our walks take nearest-neighbour steps, although this is not made explicit in the notation.  
We are interested in the velocity/speed/drift $v_1(d,\beta,\mu)$ of the first coordinate of the random walk, defined by 
\eq
\lbeq{speeddef}
v_1=\lim_{n\ra \infty}\frac{\omega_n^{[1]}}{n},
\en
whenever this limit exists.  Here $x^{[1]}$ denotes the first component of $x\in \Z^d$.

The main result of this paper is the following theorem, which is proved in Section \ref{sec-pfmainthm}.

\begin{THM}[Monotonicity of the speed]
\label{thm:main}
The velocity $v_1(d,\beta,\mu)$ for ERWD in dimension $d$ is continuous for $(\beta,\mu)\in [0,1]^2$ when $d\ge 6$ and is strictly increasing in $\beta\in [0,1]$ for each $\mu\in [0,1]$ when $d\ge 12$.
\end{THM}
In fact we show that the partial derivative of $v_1$ with respect to $\beta$ exists and is positive when $d\ge 12$.  
It is also possible to prove monotonicity in $\beta$ for $\beta +\mu$ sufficiently small when $d\ge 8$, and for all $\mu\in [0,1]$ (resp. $\beta\in [0,1]$) and $\beta$ (resp. $\mu$) sufficiently small when $d\ge 9$.  
 
In order to establish Theorem \ref{thm:main} we first need to establish that the limiting velocity $v_1$ actually exists.  Although this can be achieved for $d$ sufficiently large using the expansion of \cite{HH07}, we can appeal to the ergodic method of \cite{BSZ03}, for which it is sufficient that $d\ge 6$.  The following result is established in \cite[Theorem 1.4]{BSZ03} in the context of random walk in a random environment, however the same argument works for ERWD, and in fact much more generally (e.g. i.i.d.~random multi-cookie environment).

\begin{THM}[see Theorem 1.4 of \cite{BSZ03}]
\label{thm:speed exists}
For each $d\ge 6$ and $(\mu,\beta)\in [0,1]^2$, there exists $v(d,\mu,\beta)$ such that $\mQ(v=\lim_{n\ra \infty}n^{-1}\omega_n)=1$.
\end{THM}
In \cite{BSZ03} there exist well-behaved random cut times (independent of the environment) with the property that the set of sites visited before and after each cut-time are disjoint.  Thus one can essentially lay down a new i.i.d. environment at each of these cut-times.  The equivalent notion in the present context, would be to replace all of the cookies at these cut times since none of the sites with missing cookies at a cut time will be visited again anyway.  The following Lemma is proved in Section \ref{sec:proof-speed} using an elementary ``cookie replacement" coupling argument and implies that for $d\ge 6$ and each $\mu>0$ we can choose $\beta$ sufficiently small so that the velocity is negative.

\begin{LEM}
\label{lem:-speed}
For each $d\ge 2$ and each $\mu>0$ there exist $\delta>0$ and $\beta_*(\mu,d,\delta)>0$ such that for every $\beta<\beta_*$, $\mQ(\limsup_{n\ra \infty}\frac{\omega_n^{[1]}}{n}<-\delta)=1$.
\end{LEM}
Lemma \ref{lem:-speed} implies transience of $\omega_n^{[1]}$ for $d\ge 2$ and $\beta<\beta_*$, whence renewal techniques (e.g. \cite{BR07, SZ99,Zern05}) may be used to prove the existence of the velocity.
In high dimensions the formula for the velocity in \cite{HH07}, shows that the velocity is positive for $\beta$ sufficiently large.

\begin{LEM}
\label{lem:+speed}
For each $d\ge 9$ and each $\mu>0$ there exists $\beta^*(\mu,d)>0$ such that for every $\beta>\beta^*$, $\mQ(\lim_{n\ra \infty}\frac{\omega_n^{[1]}}{n}>0)=1$.
\end{LEM}

Lemmas \ref{lem:-speed} and \ref{lem:+speed} together with Theorem \ref{thm:main} then immediately imply the following result.
\begin{COR}
\label{cor:0speed}
For each $d\ge 9$ and each $\mu\in [0,1]$ there exists a $\beta_0(\mu,d)\in [0,1]$ for which $v=0$, i.e. $\mQ(\lim_{n\ra \infty}\frac{\omega_n}{n}=0)=1$.  For each $d\ge 12$ and $\mu\in [0,1]$ there is a unique such $\beta_0(\mu,d)$.
\end{COR}

\begin{figure}
\includegraphics[scale=.6]{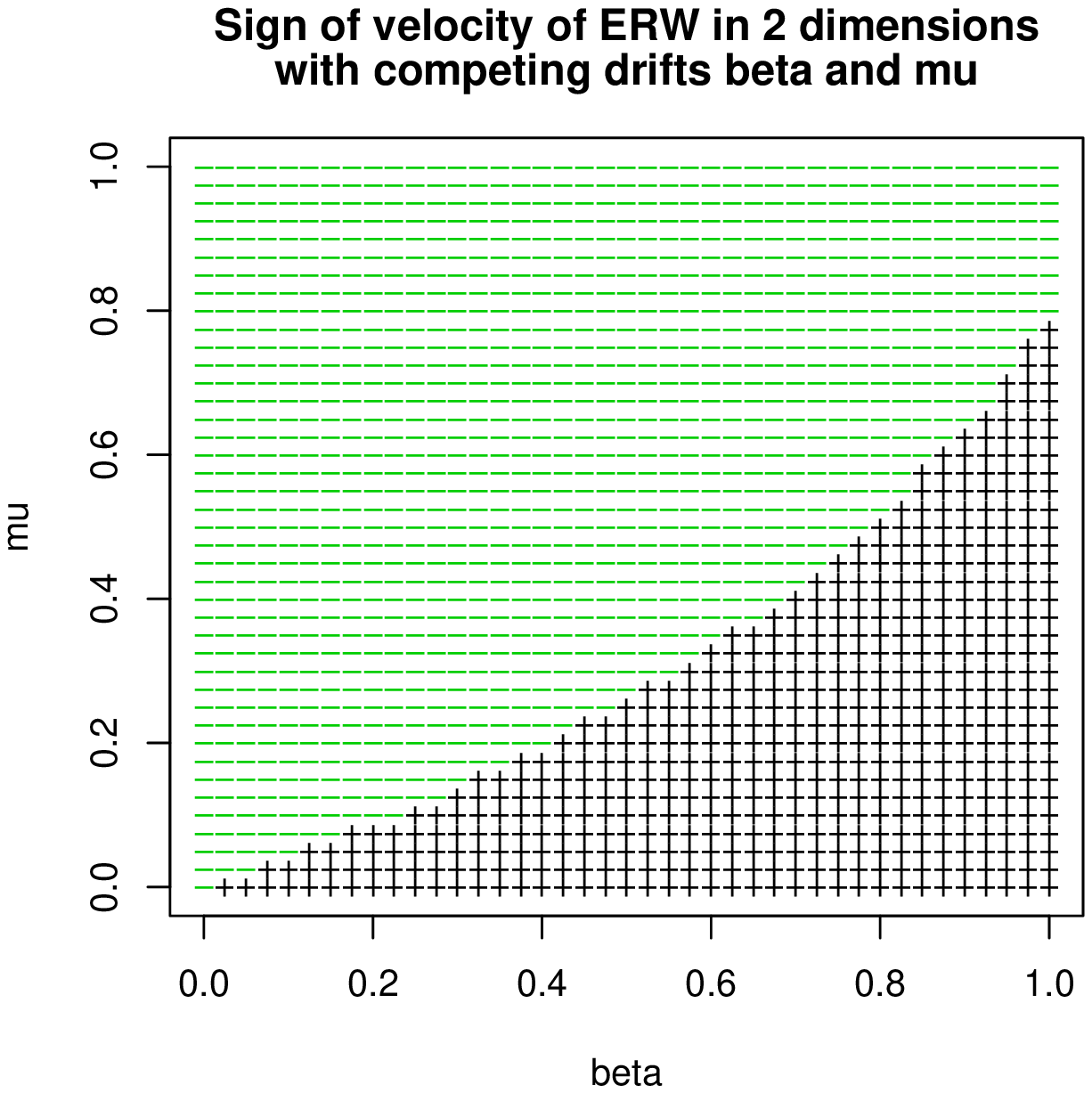}
\includegraphics[scale=.6]{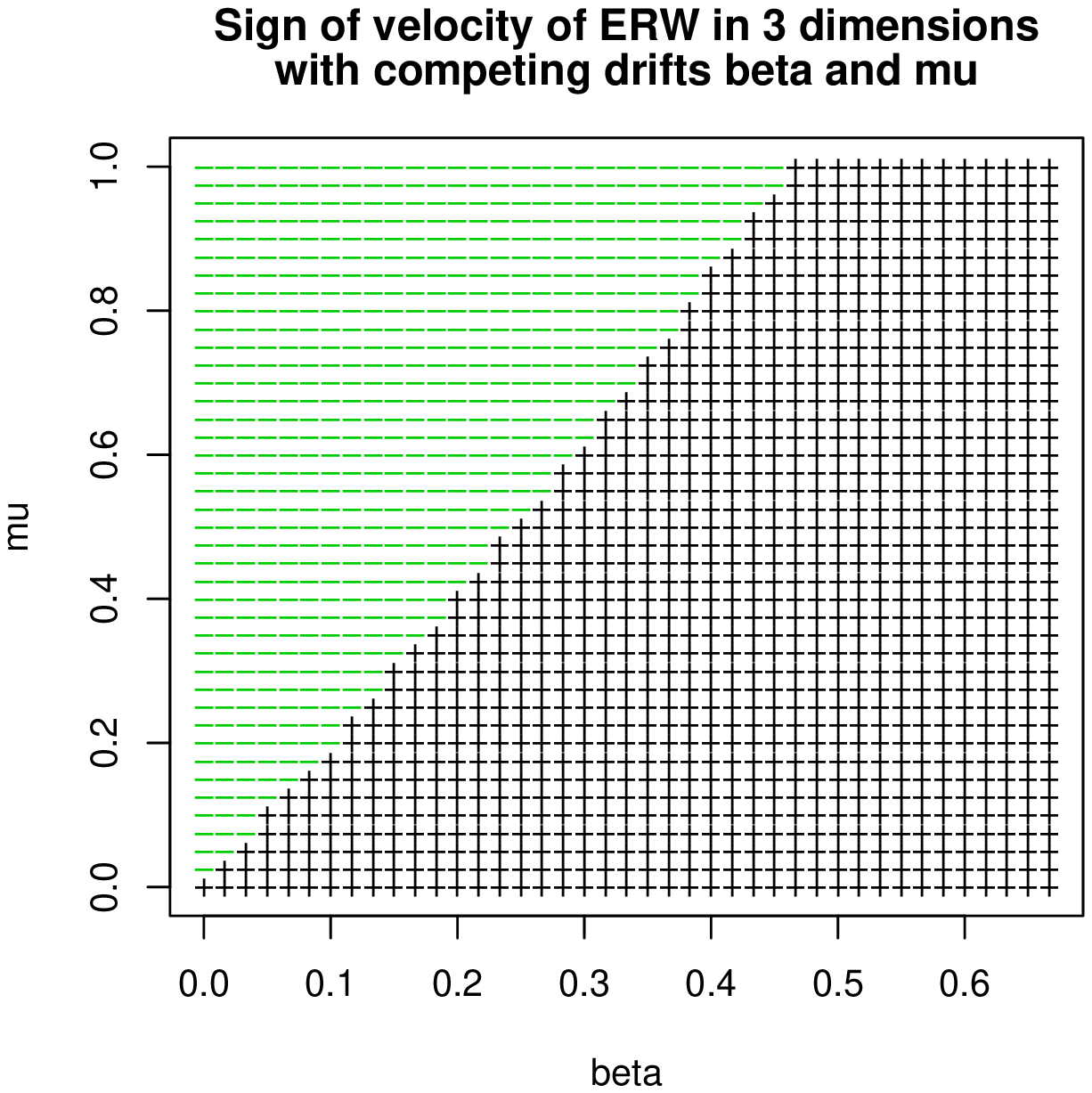}
\caption{Estimates of the sign (+,-) of the velocity $v_1$ of ERWD in 2 and 3 dimensions, based on 1000 simulations of 7000-step walks, done in R \cite{R}.} 
\label{fig}
\end{figure}

Apart from Lemma \ref{lem:-speed}, all of the above results are proved in this paper in high dimensions only.  We expect the results to hold for all $d\ge 2$, with the exception that Lemma \ref{lem:+speed} may not hold for $\mu$ close to 1, when $d=2$.  See Figure \ref{fig}. 

\blank{
d3_data=read.table("H:\\preprints\\expansion\\competing_drifts\\output3d7000_1000.txt")
pm=sign(d3_data[,3])
z=rep("-",length(pm))
z[pm>0]="+"
plot(d3_data[,1],d3_data[,2],xlab="beta",ylab="mu",main=c("Sign of velocity of ERW in 3 dimensions", "with competing drifts beta and mu"),col=-pm+2,pch=z)
}

\begin{CON}
\label{con}
For all $d\ge 2$, $(\mu,\beta)\in [0,1]^2$ the velocity $v_1=\lim_{n\ra \infty}n^{-1}\omega_n^{[1]}$ exists and is monotone increasing in $\beta$ for fixed $\mu$ and decreasing in $\mu$ for fixed $\beta$ respectively.  For each $d\ge 3$ and $\mu\in [0,1]$ there exists a unique $\beta_0(\mu,d)\in [0,1]$ such that $v(d,\mu,\beta_0)=0$.
\end{CON}
It would also be interesting to know whether the first coordinate of the walk is recurrent when $\beta=\beta_0$.  Perhaps this can be verified when $d\ge 14$ using \cite[Theorem 2.2]{BSZ03}.


The remainder of this paper is organised
as follows.  We prove Lemma \ref{lem:-speed} in Section \ref{sec:proof-speed}.  In Section \ref{sec:speed-formula} we give a formula for the speed of ERWD when 
$d\ge 6$, and verify Lemma \ref{lem:+speed}.  In Section \ref{sec-derivpi} we consider the derivative of this formula, and estimates of this derivative are used in Section \ref{sec-pfmainthm}
to prove Theorem \ref{thm:main}.

\section{Proof of Lemma \ref{lem:-speed}}
\label{sec:proof-speed}
Fix $d\ge 2$ and $\mu>0$.  We prove that for $\beta,\delta>0$ sufficiently small depending on $\mu$ and $d$, $\limsup_{n\ra \infty} n^{-1}\omega^{[1]}_{3n}<-\delta$ almost surely. 
We do this by comparing the ERWD with a similar model where all of the eaten cookies are replaced at times $3m$, $m \in\mathbb{N}$.  Perhaps the results in this section can be improved by replacing all of the cookies at times $km$ for $k>3$.

Let $E\subset \Z^d$ (resp.~$E^{-x}\subset \Z^d\setminus\{x\}$) denote a cookie configuration, i.e.~a set of vertices in $\Z^d$ (resp.~$\Z^d\setminus\{x\}$) at which cookies are present.
Let $A_x=(1-\mu+(\beta+\mu)I_{x})$ and $B_x=(1+\mu-(\beta+\mu)I_{x})$, where $I_x$ is the indicator that there is a cookie at $x$.  Then $0\le A_x,B_x\le 2$, $A_x+B_x=2$, $A_x$ is increasing in $I_x$ and $B_x$ is decreasing in $I_x$.  Let $\mQ_{E}$ denote the law of an ERWD $\vec{\eta}$ with a given initial cookie configuration $E\subset \Z^d$, and write $x\sim y$ (resp. $x\sim \sim y$) if $x$ is a neighbour of $y$ (resp. $x$ is graph distance $2$ from $y$).  Then we have
\eqalign
(2d)^3\mQ_{E}(\eta_3^{[1]}=3)=&A_0A_{e_1}A_{2e_1}\lbeq{+3}\\
(2d)^3\mQ_{E}(\eta_3^{[1]}=2)=&A_0A_{e_1}(2d-2)+A_0\sum_{u\sim e_1:u \ne e_1\pm e_1}A_u + \sum_{u \sim 0:u \ne \pm e_1}A_uA_{u+e_1}\lbeq{+2}\\
(2d)^3\mQ_{E}(\eta_3^{[1]}=1)=&
A_0A_{e_1}B_{2e_1}+
A_0B_{e_1}(1-\mu)+
B_0A_{-e_1}(1-\mu)+
A_0(2d-2)^2+\lbeq{+1_1}\\
&\sum_{u \sim 0:u\ne \pm e_1}A_u(2d-2)+
(2d-2)(1-\mu)+
\sum_{v \sim \sim 0: v\cdot e_1=0}A_v\lbeq{+1_2}\\
(2d)^3\mQ_{E}(\eta_3^{[1]}=-3)=&B_0B_{-e_1}B_{-2e_1}\lbeq{-3}\\
(2d)^3\mQ_{E}(\eta_3^{[1]}=-2)=&B_0B_{e_1}(2d-2)+B_0\sum_{u\sim -e_1:u \ne e_1\pm e_1}B_u + \sum_{u \sim 0:u \ne \pm e_1}B_uB_{u+e_1}\lbeq{-2}\\
(2d)^3\mQ_{E}(\eta_3^{[1]}=-1)=&
B_0B_{-e_1}A_{-2e_1}+
B_0A_{-e_1}(1+\mu)+
A_0B_{e_1}(1+\mu)+
B_0(2d-2)^2+\lbeq{-1_1}\\
&\sum_{u \sim 0:u\ne \pm e_1}B_u(2d-2)+
(2d-2)(1+\mu)+
\sum_{v \sim \sim 0: v\cdot e_1=0}B_v\lbeq{-1_2}\\
(2d)^3\mQ_{E}(\eta_3^{[1]}=0)=&1-\sum \text{ above terms}
\enalign
It is easy to see that \refeq{+3} and \refeq{+2} are non-decreasing in $I_x$ for every $x$, for any fixed cookie environment $E^{-x}\subset \Z^d\setminus\{x\}$, and that $\mQ_E(\eta_3^{[1]}>0)=$\refeq{+3}+\refeq{+2}+\refeq{+1_1}+\refeq{+1_2} is nondecreasing in $I_x$ for $x \notin\{0,e_1\}$ and any fixed $E^{-x}$.  Simple arithmetic enables us to write 
\eqalign
\mQ_E(\eta_3^{[1]}>0)=&A_0K_0(E^{-0})+K_0'(E^{-0})=A_{e_1}K_1(E^{-{e_1}})+K_1'(E^{-{e_1}}),
\enalign
where $K_0,K_0' \ge 0$ (resp. $K_1,K_1'\ge 0$) are constants that depend on $E^{-0}$ but not $I_0$ (resp. $E^{-e_1}$ but not $I_{e_1})$.  This verifies that \refeq{+3}+\refeq{+2}+\refeq{+1_1}+\refeq{+1_2} is nondecreasing in $I_x$ for each $x$, and each fixed $E^{-x}$.  Similarly, \refeq{-3}, \refeq{-2} and \refeq{-3}+\refeq{-2}+\refeq{-1_1}+\refeq{-1_2} are nonincreasing in $I_x$ for each $x$, $E^{-x}$.

It follows that we can define random walks $\vec{\omega}$ and $\vec{\nu}$ on the same probability space so that $\vec{\omega}$ has the law of ERWD and where $\vec{\nu}$ is the sum of independent random vectors $X_i$, each with the distribution of $\omega_3$ (corresponding to the situation above with $I_x=1$ for all $x$), and such that the following hold
\eqalign
\omega^{[1]}_{3n}-\omega^{[1]}_{3(n-1)}=3 \Rightarrow &\nu^{[1]}_{3n}-\nu^{[1]}_{3(n-1)}=3, \qquad &\nu^{[1]}_{3n}-\nu^{[1]}_{3(n-1)}=-3 \Rightarrow \omega^{[1]}_{3n}-\omega^{[1]}_{3(n-1)}=-3,\nn\\
\omega^{[1]}_{3n}-\omega^{[1]}_{3(n-1)}=2 \Rightarrow &\nu^{[1]}_{3n}-\nu^{[1]}_{3(n-1)}\ge 2, \qquad &\nu^{[1]}_{3n}-\nu^{[1]}_{3(n-1)}=-2 \Rightarrow \omega^{[1]}_{3n}-\omega^{[1]}_{3(n-1)}\le -2,\nn\\
\omega^{[1]}_{3n}-\omega^{[1]}_{3(n-1)}=1 \Rightarrow &\nu^{[1]}_{3n}-\nu^{[1]}_{3(n-1)}\ge 1, \qquad &
\nu^{[1]}_{3n}-\nu^{[1]}_{3(n-1)}=-1 \Rightarrow \omega^{[1]}_{3n}-\omega^{[1]}_{3(n-1)}\le -1\nn.
\enalign
In particular, under this coupling, $\nu_{3n}^{[1]}-\omega_{3n}^{[1]}$ is non-negative and non-decreasing in $n$, almost surely.

Since $n^{-1}\sum_{i=1}^n X_i\ra \mE[\omega_3]$ almost surely, and $\mE[\omega_3^{[1]}]$ is continuous in $\beta$, it is enough to prove that $\mE[\omega_3^{[1]}]<0$ when $\beta=0$.  In this case $A_x=B_x=1$ for each $x$ so we can use \refeq{+3}-\refeq{-1_2} to get $(2d)^{3}\mE[\omega_3^{[1]}]=-4d\mu<0$ when $\beta=0$ as required.  \qed

Note that by keeping track of the $\beta$ dependence in \refeq{+3}-\refeq{-1_2}, one can get an explicit bound on how small $\beta$ should be as a function of $\mu$ and $d$.

\section{The formula for the speed}
\label{sec:speed-formula}
We recall some notation and results from \cite{HH07} and \cite{HH08mono} giving an expression for the velocity $v$ of the ERWD.
If $\ewalk$ and $\walk$ are two paths of length at least $j$ and $m$
respectively and such that $\eta_j=\omega_0$, then the concatenation
$\ewalk_j\circ \walk_m$ is defined by
    \eq
    \lbeq{concat}
    (\ewalk_j\circ \walk_m)_{\ch{i}}=\left\{
    \begin{array}{lll}
    &\eta_i&\text{when }0\le i\leq j,\\
    &\omega_{i-j} &\text{when }j \leq i \leq m+j.
    \end{array}\right.
    \en
Given $\ewalk_m$ such that $\Qbold(\walk_m=\ewalk_m)>0$, we define a conditional probability measure $\Qbold^{\ewalk_m}$
on walks starting from $\eta_m$ by 
    \eq
    \lbeq{cylinders}
    \Qbold^{\ewalk_m} (\walk_n=\vec{x}_n)
    \equiv\prod_{i=0}^{n-1} p^{\ewalk_m\circ \ci{\vec{x}_i}}(x_{i},x_{i+1})=\Qbold(\walk_{m+n}=\ewalk_m\circ \vec{x}_n|\walk_m=\ewalk_m).
    \en

Set $j_0=0$, and for $N\geq 1$ let
    \eq
    \lbeq{DeltaNdef}
    \Delta_{\sss N}= \big(p^{\walkvec{N-1}{j_{N-1}+1}\circ
    \walkvec{N}{j_{N}}}-p^{\walkvec{N}{j_{N}}}\big)(\walkcoor{N}{j_{N}},\walkcoor{N}{j_{N}+1}).
    \en
For ERWD, $\Delta_{\sss N}$ is non-zero precisely when $\walkcoor{N}{j_{\sN}}$ has already been visited by $\walkvec{N-1}{j_{N-1}+1}$ but not by $\walkvec{N}{j_{N}-1}$, so that
    \eqalign
    \lbeq{Deltabound}
    \Delta_{\sss N}=&\left|\frac{(\beta+\mu) e_1\cdot(\walkcoor{N}{j_{N}+1}-\walkcoor{N}{j_{N}})}{2d}\left[I_{\{\walkcoor{N}{j_{N}}\notin \walkvec{N-1}{j_{N-1}}\circ \walkvec{N}{j_{\sN}-1}\}}-I_{\{\walkcoor{N}{j_{N}}\notin \walkvec{N}{j_{\sN}-1}\}}\right]\right|, \text{ and}\nn\\
|\Delta_{\sss N}|    \le &
    \frac{\beta+\mu}{2d}I_{\{\walkcoor{N}{j_{N}+1}=\walkcoor{N}{j_{N}}\pm e_1\}}I_{\{\walkcoor{N}{j_{N}}\in \walkvec{N-1}{j_{N-1}}\}}.
    \enalign

Define $\mc{A}_{m,\sN}=\{(j_1, \dots, j_{\sN})\in \Z_+^{\sN}:\sum_{l=1}^Nj_l=m-N-1\}$, $\mc{A}_{\sN}=\bigcupd_{m}\mc{A}_{m,\sN}$ and
      \eqalign
    \pi_m^{\smallsup{N}}(x,y):=&\sum_{\vec{j}\in \mc{A}_{m,N}}\sum_{\walkvec{0}{1}}\sum_{\walkvec{1}{j_{\sss 1}+1}}\dots\sum_{\walkvec{N}{j_{\sss N}+1} }I_{\{\omega^{(N)}_{j_{\sN}}=x, \omega^{(N)}_{j_{\sN}+1}=y\}}p^0(0,\walkcoor{0}{1})
    \prod_{n=1}^{N}\Delta_{\sss n}\prod_{i_{n}=0}^{j_{n}-1}p^{\walkvec{n-1}{j_{n-1}+1}\circ \walkvec{n}{i_{n}}}\left(\walkcoor{n}{i_{n}},\walkcoor{n}{i_{n}+1}\right).
    \lbeq{piNxydef}
    \enalign
Note that this quantity depends on $\beta$, $\mu$, and $d$, as do the following
    \eq
    \pi_m(x,y):=\sum_{N=1}^{\infty} \pi_m^{\smallsup{N}}(x,y),  \quad
    \pi^{\smallsup{N}}(x,y):=\sum_{m}\pi_m^{\smallsup{N}}(x,y),
    \quad  \text{and} \quad \pi_m(y):=\sum_{N=1}^{\infty}\sum_x \pi_m^{\smallsup{N}}(x,y).
    \lbeq{piotherdef}
    \en
In \cite{HH07}, it was shown that if $\lim_{n \ra \infty}\sum_{m=2}^n \sum_x x\pi_m(x)$ exists
and $n^{-1}\omega_n\ra v$ in probability, then
\eqalign
    v(\beta,\mu,d)& = \sum_x xp^0(0,x) +\sum_{m=2}^\infty \sum_x x\pi_m(x)\lbeq{theta}\\
    & =\frac{\beta e_1}{d}+\sum_{m=2}^\infty\sum_{N=1}^{\infty}\sum_{x,y}  (y-x)\pi_m^{\smallsup{N}}(x,y)\lbeq{theta2}.
    \enalign




Let $\mP_{d}$ denote the law of simple symmetric random walk in
$d$ dimensions, beginning at the origin, and let $D_d(x)=(2d)^{-1}I_{\{|x|=1\}}$ be the
simple random walk step distribution.  For absolutely summable functions $f,g$ on $\Zd$, define the convolution of $f$ and $g$ by  
\eq
    (f*g)(x) = \sum_y f(y) g(x-y).
    \en
Let $f^{*k}(x)$ denote the
$k$-fold convolution of $f$ with itself, and let
$G_d(x)=\sum_{k=0}^{\infty} D_d^{*k}(x)$ denote the Green's
function for this random walk.   Then $G_d^{*k}(x)<\infty$ when $d>2k$.  Define $G_d^{*k}=G_d^{*k}(0)$.

For $i\geq 0$ and $d>2(i+1)+1$ define
    \eqalign
    \mc{E}_i(d)=&
    \big(\frac{d}{d-1}\big)^{i+1}G_{d-1}^{*(i+1)}-1.\lbeq{Eidef}
    \enalign
For $d>5$ define
    \eqalign
    a_d=\frac{d}{(d-1)^{2}}G^{*2}_{d-1}\lbeq{addef}.
    \enalign

The following proposition is proved exactly as in the proof of \cite[Proposition 3.2]{HH08mono} with $\beta$ replaced by $\beta+\mu$ in each of the bounds.
\begin{PRP}[Bounds on the expansion coefficients]
\label{prp:pibound}
For all $N\ge 1$,
\eqalign
\sum_{x,y}\sum_m|\pi_m^{\smallsup{N}}(x,y)|\le
\begin{cases}
(\beta+\mu) d^{-1}\mc{E}_0(d) & N=1,\\
(\beta+\mu)^Nd^{-1}(d-1)^{-1}G_{d-1}\mc{E}_1(d)\left(a_d\right)^{N-2} & N>1.\lbeq{pibound}
\end{cases}
\enalign
\end{PRP}

Since the speed is known to exist when $d\ge 6$ by the ergodic argument of \cite{BSZ03} (see Theorem \ref{thm:speed exists}), the following corollary is an easy consequence of
\cite[Propositions 3.1 and 6.1]{HH07}
together with Proposition \ref{prp:pibound}, and the fact that $2a_6<1$ since $G_5^{*2}<25/12$ \cite{HS92b}.
\begin{COR}
\label{cor:ERWspeed}
For all $d\ge 6$ and $\beta,\mu \in [0,1]$,
    \eqalign
    v(\beta,\mu,d)
    =\frac{\beta e_1}{d}+\sum_{m=2}^{\infty}\sum_x x \pi_m(x)=\frac{\beta e_1}{d}+\sum_{N=1}^{\infty}\sum_{m=2}^{\infty}\sum_{x,y} (y-x) \pi^{(N)}_m(x,y).
    \enalign
\end{COR}

\noindent {\em Proof of Lemma \ref{lem:+speed}.}  Corollary  \ref{cor:ERWspeed} implies that for $d\ge 6$, 
\eq
\big|v_1-\frac{\beta}{d}\big|\le \sum_{N=1}^{\infty}\sum_{x,y}\sum_m|\pi_m^{\smallsup{N}}(x,y)|.
\en
Thus to prove the lemma, it is sufficient to show that there exists $\delta>0$ independent of $\beta$ such that $d\sum_{N=1}^{\infty}\sum_{x,y}\sum_m|\pi_m^{\smallsup{N}}(x,y)|<1-\delta$.  From Proposition \ref{prp:pibound} we have 
\eq
d\sum_{N=1}^{\infty}\sum_{x,y}\sum_m|\pi_m^{\smallsup{N}}(x,y)|\le 2\mc{E}_0(d)+4\frac{G_{d-1}\mc{E}_1(d)}{(d-1)(1-2a_d)},
\en
and since the right hand side is independent of $\mu$ and $\beta$ and decreasing in $d$, it is enough to prove that it is less than 1 when $d=9$.  This result can easily be checked using the rigorous upper bounds $G_8\le 1.07865$ and $G_8^{*2}\le 1.28901$ \cite{Hpc,HS92b}.

\blank{G1vec7plus<-c(1.09390631558785,1.07864701201693,1.06774608638140,1.05954374788826,1.05313615290863)
G2vec7plus<-c(1.36678616952622,1.28900278970224,1.23987303101441,1.20560748224037,1.18017394904878)
G3vec7plus<-c(2.35756693371046,1.83154611054153,1.62536150147784,1.50817443085253,1.43042772680824)
#D2G2vec7plus<-c(.179,.132,.105,.087,.074)
#supG2minusdelta<-pmax(G2vec7plus-1,D2G2vec7plus)
#G1neighbourupperbound<-rep(.0051,5)

dtest<-function(d){
G1<-G1vec7plus[d-7]
G2<-G2vec7plus[d-7]
G3<-G3vec7plus[d-7]
E_0<-((d/(d-1))*G1)-1
E_1<-((d/(d-1))^2*G2)-1
eps<-2*d*G1*G3/((d-1)^4)+E_1*G2/(d*(d-1)^2)
a_d<-d*G2/((d-1)^2)
#drho<-E_0/d+G1*E_1/(d*(d-1)*(1-a_d))
#dchi<-E_0 + G1*E_1*(2-a_d)/((d-1)*((1-a_d)^2))
#dgamma<-d*G2/(d-1)^2+(eps*d)/(1-a_d)+2*d*E_1*G1*G3/((d-1)^4 * (1-a_d)^2)
val=2*E_0+4*G1*E_1/((d-1)*(1-2*a_d))
#val<-drho+dgamma+dchi
list("want this <1"=val,"E_0"=E_0)}
}

\qed

\section{The differentiation step}
\label{sec-derivpi}
We follow the analysis in \cite{HH08mono}.
We fix $\mu\in [0,1]$, differentiate
the right hand side of \refeq{theta2} with respect to $\beta$, and prove that this derivative
is positive for all $\beta\in [0,1]$, when $d\ge 12$.
Letting $\varphi_m^{\smallsup{N}}(x,y)=\frac{\del}{\del\beta}\pi_m^{\smallsup{N}}(x,y)$
and assuming that the limit can be taken through the infinite sums, we then have
    \eqalign
    \lbeq{deriv1}
    \frac{\del v_1}{\del \beta}(\beta,\mu,d)&=\frac{1}{d}+\sum_{N=1}^{\infty}\sum_{m=2}^\infty\sum_{x,y}  (y_1-x_1)\varphi_m^{\smallsup{N}}(x,y).
    \enalign
Since $\varphi_m^{\smallsup{N}}(x,y)\equiv 0$ unless $|x-y|=1$, we have that
    \eqalign
    \left|\frac{\del v_1}{\del \beta}(\beta,d)-\frac{1}{d}\right|\le \sum_{N=1}^{\infty}\sum_{m=2}^\infty\sum_{x,y}|\varphi_m^{\smallsup{N}}(x,y)|.\lbeq{needed1}
    \enalign
We conclude that $\frac{\del v_{1}}{\del \beta}$ 
is positive when $\sum_{N=1}^{\infty}\sum_{m=2}^\infty\sum_{x,y}|\varphi_m^{\smallsup{N}}(x,y)|<d^{-1}$.
To verify the exchange of limits in \refeq{deriv1},
it is sufficient to prove that $\sum_{x,y}(y-x)\pi_m^{\smallsup{N}}(x,y)$
is absolutely summable in $m$ and $N$ and that
$\sum_{N=1}^{\infty}\sum_{m=2}^\infty\sup_{\beta,\mu\in [0,1]}|\sum_{x,y}  (y-x)\varphi_m^{\smallsup{N}}(x,y)|<\infty$.
By Proposition \ref{prp:pibound} and the fact that $|y-x|=1$
for $x,y$ nearest neighbours, the first condition holds provided that
    \eqalign
    \boxed{(\beta+\mu) a_d<1.}\lbeq{C1}
    \enalign
In fact we shall see that this inequality for $\beta=\mu=1$ is also sufficient to also establish the second condition. We now identify
$\varphi_m^{\smallsup{N}}(x,y)$.

It follows from \refeq{piNxydef} that
    \eqalign
    \varphi_m^{\smallsup{N}}(x,y)=\varphi_m^{\smallsup{N,1}}(x,y)+\varphi_m^{\smallsup{N,2}}(x,y)+\varphi_m^{\smallsup{N,3}}(x,y),\lbeq{varphibreak}
    \enalign
where (by Leibniz' rule), $\varphi_m^{\smallsup{N,1}}(x,y)$, $\varphi_m^{\smallsup{N,2}}(x,y)$ and $\varphi_m^{\smallsup{N,3}}(x,y)$ arise from differentiating $p^0(0,\walkcoor{0}{1})$, $\prod_{n=1}^{N}\prod_{i_{n}=0}^{j_{n}-1}p^{\walkvec{n-1}{j_{n-1}+1}\circ \walkvec{n}{i_{n}}}\left(\walkcoor{n}{i_{n}},\walkcoor{n}{i_{n}+1}\right)$ and $\prod_{n=1}^{N}\prod_{i_{n}=0}^{j_{n}-1}\Delta_{\sss N}$, respectively.

Observe that if $\vec{\eta}_{m}=x_l$ then
    \eqalign
    \frac{\del}{\del\beta}{p_{\beta}^{\vec{\eta}_m}}(x_{l},x)=&\frac{I_{\{x_l \notin\vec{\eta}_{m-1}\}}}{2d}\left(I_{\{x-x_l= e_1\}}-I_{\{x-x_l= -e_1\}}\right)\lbeq{pderiv},
    \enalign
and hence, using $I_{A}-I_{A\cap C}=I_{A \cap C^c}$ we have
    \eqalign
    \frac{\del}{\del\beta}\left({p_{\beta}^{\vec{\eta}_m}}(x_{l},x)-{p_{\beta}^{\vec{\omega}_n\circ\vec{\eta}_m}}(x_{l},x)\right)
    =\frac{1}{2d}I_{\{x_l \notin\vec{\eta}_{m-1},x_l \in \vec{\omega}_{n-1}\}}\left(I_{\{x-x_l= e_1\}}-I_{\{x-x_l= -e_1\}}\right).\lbeq{deltaderiv}
    \enalign
Clearly then
    \eqalign
    \left|\frac{\del}{\del\beta}\left({p_{\beta}^{\vec{\eta}_m}}(x_{l},x)-{p_{\beta}^{\vec{\omega}_n\circ\vec{\eta}_m}}(x_{l},x)\right)\right|
    \le \frac{1}{2d}I_{\{x_l \in \vec{\omega}_{n-1}\setminus \vec{\eta}_{m-1}\}}\left(I_{\{x-x_l= e_1\}}+I_{\{x-x_l= -e_1\}}\right).\lbeq{deltaderivrep}
    \enalign

Let $\rho^{\smallsup{N}}$ be obtained by replacing
$p^0(0,\walkcoor{0}{1})$ in \refeq{piNxydef} with $(2d)^{-1}I_{\{\walkcoor{0}{1}= \pm e_1\}}$ (a bound on its derivative)
 and by bounding $\Delta_{\sss N}$ by $|\Delta_{\sss N}|$
for all $n=1, \ldots, N$.

For $k=1, \ldots, N$, let $\gamma^{\smallsup{N}}_k$
be obtained from \refeq{piNxydef} by bounding
$\Delta_{\sss N}$ by $|\Delta_{\sss N}|$
for all $n=1, \ldots, N$ and by replacing $\prod_{i_{k}=0}^{j_{k}-1}p^{\walkvec{k-1}{j_{k-1}+1}\circ \walkvec{k}{i_{k}}}\left(\walkcoor{k}{i_{k}},\walkcoor{k}{i_{k}+1}\right)$ with the following bound on its derivative
\eq
\sum_{l=0}^{j_{k}-1}\frac{I_{\{\walkcoor{k}{l_{k}+1}-\walkcoor{k}{l_{k}}= \pm e_1\}}}{2d}\prodtwo{i_{k}=0}{i_{k}\ne l}^{j_{k}-1}p^{\walkvec{k-1}{j_{k-1}+1}\circ \walkvec{k}{i_{k}}}\left(\walkcoor{k}{i_{k}},\walkcoor{k}{i_{k}+1}\right).
\en
Similarly, let $\chi^{\smallsup{N}}_k$
be obtained by replacing $\Delta^{\smallsup{k}}_{j_{k}+1}$ in
\refeq{piNxydef} by
$(2d)^{-1}I_{\{\walkcoor{k}{j_{k}} \in \walkvec{k-1}{j_{k-1}+1}\setminus\walkvec{k}{j_{k}-1}\}}
I_{\{\walkcoor{k}{j_{k}+1}-\walkcoor{k}{j_{k}}= \pm e_1\}}$ (a bound on its derivative) and by bounding
$\Delta_{\sss n}$ for $n\neq k$ by $|\Delta_{\sss n}|$.

Letting $\gamma^{\smallsup{N}}=\sum_{k=1}^{N}\gamma^{\smallsup{N}}_k$ and
$\chi^{\smallsup{N}}=\sum_{k=1}^{N}\chi^{\smallsup{N}}_k$, we obtain that
    \eqalign
    \sum_{m}\sum_{x,y}|\varphi_m^{\smallsup{N,1}}(x,y)|\le \rho^{\smallsup{N}}, \quad \sum_{m}\sum_{x,y}|\varphi_m^{\smallsup{N,2}}(x,y)|\le
    \gamma^{\smallsup{N}}, \quad \text{and  }\sum_{m}\sum_{x,y}|\varphi_m^{\smallsup{N,3}}(x,y)|\le
    \chi^{\smallsup{N}}.\lbeq{3terms}
    \enalign
Define
    \eqn{
    \epsilon(d)=
    \frac{2d}{(d-1)^4}G_{d-1}G_{d-1}^{*3}+\frac{\mc{E}_1(d)}{d(d-1)^2}G_{d-1}^{*2},
    \lbeq{epsilonddef}
    }
The following proposition is proved exactly as in Corollary 4.5 of \cite{HH08mono} except that we replace $\beta$ with $\beta+\mu$ in all of the bounds.
\begin{PRP}[Summary of bounds]
\label{prp:derivbounds}
For all $\beta,\mu \in [0,1]$, and $d$ such that $2a_d<1$,
    \eqalign
    d\sum_{N=1}^{\infty}\rho^{(N)}\le
    &\frac{2\mc{E}_0(d)}{d}+\frac{4G_{d-1}\mc{E}_1(d)}{d(d-1)\big(1-2a_d\big)}\lbeq{RHO}\\
    d\sum_{N=1}^{\infty}\chi^{(N)}\le &\mc{E}_0(d)+\frac{2G_{d-1}\mc{E}_1(d)(2-2a_d)}{(d-1)\big(1-2a_d\big)^2}\lbeq{CHI}\\
    d\sum_{N=1}^{\infty}\gamma^{(N)}\le & \frac{2dG_{d-1}^{*2}}{(d-1)^2}+\frac{4\epsilon(d)d}{\big(1-2a_d\big)}+
    \frac{16d\mc{E}_1(d)G_{d-1}G_{d-1}^{*3}}{(d-1)^4\big(1-2a_d\big)^{2}}.\lbeq{GAMMA}
    \enalign
    
\end{PRP}

\section{Proof of Theorem \ref{thm:main}}
\label{sec-pfmainthm}
Continuity of $v_1$ as a function of $\mu$ and $\beta$ follows from the fact that the formula for the speed is a sum of functions that are continuous in $\mu$ and $\beta$ and absolutely summable, uniformally in $(\mu,\beta)\in [0,1]^2$ by Proposition \ref{prp:pibound} and the fact that $2a_d<1$ for $d\ge 6$. 

For $d\ge 6$, the bounds of Proposition \ref{prp:derivbounds}
hold.  From \refeq{3terms} we have the required absolute summability
conditions in the discussion after \refeq{deriv1}, and in particular
\refeq{deriv1} holds for all $\beta$.  To complete the proof of the
theorem, it remains to show that the right hand side of \refeq{needed1}
is no more than $d^{-1}$.  By \refeq{3terms} and Proposition \ref{prp:derivbounds},
we have bounded $d$ times the right hand side of \refeq{needed1} by the sum
of the right hand sides of the bounds in Proposition \ref{prp:derivbounds}.
Since these terms all involve simple random walk Green's functions
quantities, we will need to use estimates of these quantities.


By \cite[Lemma C.1]{HS92b}, $d\mapsto G_{d}^{*n}$
is monotone decreasing in $d$ for each $n\geq 1$, so that it
suffices to show that the sum of terms on the
right hand sides of \refeq{RHO}, \refeq{CHI} and \refeq{GAMMA}
is bounded by $1$ for $d=12$.
For this we use the following rigorous Green's
functions upper bounds \cite{Hpc, HS92b}:
    \ci{\eqalign
    &G_{11}(0)\le 1.05314, \quad G^{*2}_{11}(0)\le 1.18018, \quad G^{*3}_{11}(0)\le 1.43043.
    \enalign
    }
    
\blank{G1vec7plus<-c(1.09390631558785,1.07864701201693,1.06774608638140,1.05954374788826,1.05313615290863)
G2vec7plus<-c(1.36678616952622,1.28900278970224,1.23987303101441,1.20560748224037,1.18017394904878)
G3vec7plus<-c(2.35756693371046,1.83154611054153,1.62536150147784,1.50817443085253,1.43042772680824)
#D2G2vec7plus<-c(.179,.132,.105,.087,.074)
#supG2minusdelta<-pmax(G2vec7plus-1,D2G2vec7plus)
#G1neighbourupperbound<-rep(.0051,5)

dtest<-function(d){
G1<-G1vec7plus[d-7]
G2<-G2vec7plus[d-7]
G3<-G3vec7plus[d-7]
E_0<-((d/(d-1))*G1)-1
E_1<-((d/(d-1))^2*G2)-1
eps<-2*d*G1*G3/((d-1)^4)+E_1*G2/(d*(d-1)^2)
a_d<-d*G2/((d-1)^2)
drho<-2*E_0/d+4*G1*E_1/(d*(d-1)*(1-2*a_d))
dchi<-E_0 + 2*G1*E_1*(2-2*a_d)/((d-1)*((1-2*a_d)^2))
dgamma<-2*d*G2/(d-1)^2+4*(eps*d)/(1-2*a_d)+16*d*E_1*G1*G3/((d-1)^4 * (1-2*a_d)^2)
val<-drho+dgamma+dchi
list("want this <1"=val,"E_0"=E_0)}
}    
    
Putting in these values we get that the sum of the
right hand sides of the bounds in Proposition \ref{prp:derivbounds}
is at most $0.847$, for $d\geq 12$.

\qed

\paragraph{Acknowledgements.}
This work was supported in part by a FRDF grant from the University of Auckland.  We would like to thank Takashi Hara for providing the SRW Green's functions upper bounds.


\end{document}